\documentclass[10pt]{article}

\usepackage{a4wide}
\usepackage{amssymb}
\usepackage{amsfonts}
\usepackage{amsmath}
\input xy
\xyoption{arrow} \xyoption{matrix}

\date{}

\newtheorem{proposition}{Proposition}[section]
\newtheorem{theorem}[proposition]{Theorem}
\newtheorem{lemma}[proposition]{Lemma}

\def\der{\partial }

\def\nFM0{{\nu }_{F,M_0}}
\def\nFN0{{\nu }_{F,N_0}}
\def\nGN0{{\nu }_{G,N_0}}

\def\N0{ {\bf N}_0 }

\def\t{\otimes}

\def\v{\varphi}
\def\ra{\rightarrow}

\def\lra{\leftrightarrow}
\def\Xpm{X^{\pm }}

\def\s{\sigma}
\def\Z{\mathbb{Z}}

\def\l1{{\lambda}_1}

\def\a{\alpha}
\def\a0{ {\alpha }_0}
\def\a1{ {\alpha }_1}

\def\l{\lambda}
\def\o{\omega}

\def\nFGM0{{\nu }_{F,G,M_0}}

\def\nFN0{{\nu}_{F,N_0}}

\def\sm{{\sigma}^m}

\def\sm1{{\sigma}^{-1}}

\def\smtp1{{\sigma}^{-t+1}}

\def\o{\omega }
\def\S1{S^{-1}}

\def\Xpm1{X^{\pm 1}_1}

\def\sPM1{{\sigma }^{\pm 1}}
\def\sMP1{{\sigma }^{\mp 1 }}

\def\d{\delta}

\def\di{{\rm d.ind}}

\def\L{\Lambda}

\def\Ytm1{Y^{t-1}}
\def\Yim1{Y^{i-1}}

\def\CM{{\cal M}}

\def\Aut{{\rm Aut}}

\def\ad{{\rm ad }}
\def\dim{{\rm dim }}

\def\ker{ {\rm ker } }

\def\Ev{ {\rm Ev} }

\def\SL2Z{ {\rm SL}_2({\bf Z}) }

\def\th{ \theta }

\def\Gp1{ G^{1 , 1 } }
\def\P11{ P^{-1 , 1 } }
\def\Pp1{ P^{1 , 1 } }

\def\th{\theta}

\def\nCLsr{{}^\nu\kern-2pt {\cal L}^{\sigma , \rho  }}
\def\nP{{}^\nu \kern-2pt P}
\def\nL{{}^\nu\kern-2pt L}
\def\nLL{{}^\nu\kern-2pt \Lambda}
\def\nPsr{{}^\nu\kern-2pt P^{\sigma , \rho  }}
\def\nLsr{{}^\nu\kern-2pt L^{\sigma , \rho  }}
\def\nuCL{{}^\nu\kern-2pt  {\cal L}}
\def\nCLsr{{}^\nu\kern-2pt {\cal L}^{\sigma , \rho  }}
\def\nCL1m{{}^\nu\kern-2pt {\cal L}^{-1 , 1  }}
\def\x1nu{x^\frac{1}{\nu}}
\def\xm1nu{x^{-\frac{1}{\nu}}}

\def\ra{\rightarrow }

\def\CB{{\cal B}}

\def\coker{{\rm coker}}

\def\nAM0{{\nu }_{{\cal A},M_0}}
\def\nAN0{{\nu }_{{\cal A},N_0}}

\def\End{ {\rm End }}

\def\det{ {\rm det }}
\def\ad{ {\rm ad }}

\def\GL{{\rm GL}}
\def\SL{{\rm SL}}

\def\di!{\frac{\der^i}{i!}}
\def\dik!{\frac{\der^k_i}{k!}}

\def\N{\mathbb{N}}

\def\0{\overline{0}}
\def\1{\overline{1}}

\def\Ln1{\L_{n,\overline{1}}}

\def\oa{\overline{a}}

\def\a1{a_{\overline{1}}}

\def\S{\Sigma}

\def\vn1{\overrightarrow{n-1}}

\def\im{{\rm im}}

\def\Inn{{\rm Inn}}

\def\mJ{\mathbb{J}}
\def\mI{\mathbb{I}}

\def\Cen{{\rm Cen}}

\def\mT{\mathbb{T}}
\def\ind{{\rm ind}}

\def\K1{{\rm K}_1}

\def\rG{\mathrm{G}}

\def\hmI1{\widehat{\mI_1}}
\def\tmI1{\widetilde{\mI_1}}
\def\tmJ1{\widetilde{\mJ_1}}
\def\hB1{\widehat{B_1}}
\def\hCB1{\widehat{\CB_1}}

\begin{document}

\author{V. V. \  Bavula 
}

\title{An analogue of the Conjecture of Dixmier is true for the algebra
 of polynomial integro-differential operators}

\maketitle

\begin{abstract}
Let  $A_1:=K\langle x, \frac{d}{dx} \rangle$ be the Weyl algebra
and $\mI_1:= K\langle x, \frac{d}{dx}, \int \rangle$ be the
algebra of polynomial integro-differential operators over a field
$K$ of characteristic zero. The Conjecture/Problem of Dixmier
(1968) [still open]: {\em is an algebra endomorphism of the Weyl
algebra $A_1$ an automorphism?} The aim of the paper is to prove
that {\em each algebra endomorphism of the algebra $\mI_1$ is an
automorphism}. Notice that in contrast to the Weyl algebra $A_1$
the algebra $\mI_1$ is a non-simple, non-Noetherian  algebra which
is not a domain. Moreover, it contains infinite direct sums of
nonzero left and right  ideals.



\end{abstract}


\section{Introduction}

In this paper, $A_1:=K\langle x, \frac{d}{dx} \rangle$ is the {\em
Weyl algebra} (i.e. the algebra of polynomial differential
operators)  and $\mI_1:= K\langle x, \frac{d}{dx}, \int \rangle$
is the {\em algebra of polynomial integro-differential operators}
over a field $K$ of characteristic zero ($A_1, \mI_1\subseteq
\End_K(K[x])$ where $K[x]$ is a polynomial algebra in one variable
$x$), and $\int :K[x]\ra K[x]$, $x^n\mapsto \frac{x^{n+1}}{n+1}$,
$n\geq 0$, is the {\em integration}.

$\noindent $

{\bf Six Problems of Dixmier, \cite{Dix}, for the Weyl algebra
$A_1$}: In 1968, Dixmier \cite{Dix} posed six problems for the
Weyl algebra $A_1$.

$\noindent $

{\bf The First Problem/Conjecture of Dixmier, \cite{Dix}}: {\em is
an algebra endomorphism of the Weyl algebra $A_1$ an
automorphism?}

$\noindent $

Dixmier writes in his paper \cite{Dix}, p. 242: ``{\em A. A.
Kirillov informed me that the Moscow school also considered this
problem}".

$\noindent $

In 1975, the Third Problem of Dixmier was solved by Joseph and
Stein \cite{josclA1} (using results of McConnel and Robson
\cite{MR-JA-1973}); and  using his (difficult) polarization
theorem for the  Weyl algebra $A_1$  Joseph \cite{josclA1} solved
the Sixth Problem of Dixmier (a short proof to this problem is
given in \cite{BavDP6}; moreover, an analogue of the Sixth Problem
of Dixmier is true  for the ring of differential operators on an
arbitrary smooth irreducible algebraic curve \cite{BavDP6}). In
2005, the Fifth Problem of Dixmier was solved in \cite{BavDP5}.
Problems 1, 2, and 4 are still open. The Fourth Problem of Dixmier
has positive solution for {\em all homogeneous} elements of the
Weyl algebra $A_1$ (Theorem 2.3, \cite{BavDP5}).

$\noindent $

The aim of the paper is to prove an analogue of the First
Problem/Conjecture of Dixmier for the algebra $\mI_1$ (Theorem
\ref{11Oct10}). The proof is not straightforward and several key
results of the  papers \cite{Bav-algintdif}, \cite{Bav-indifaut}
and \cite{Bav-intdifline} are used. To make the proof more
accessible for the reader
 we use a `zoom in' way of presenting it: in the Introduction  we
 explain the structure of the proof, it consists of nine steps;  in
 Section \ref{PTH11} each steps is proved using some of the results
  of  \cite{Bav-algintdif}, \cite{Bav-indifaut} and \cite{Bav-intdifline}.
\begin{theorem}\label{11Oct10}
Each algebra endomorphism of $\mI_1$ is an automorphism.
\end{theorem}

{\it Structure of the Proof}. Let $\s$  be an algebra endomorphism
of $\mI_1$. Since $\mI_1=K \langle H, \int , \der \rangle $ where
$H:=\der x$ (notice that $x=\int H$), the endomorphism $\s$ is
uniquely determined by the elements
$$H':=\s (H), \;\; \int':= \s (\int ), \;\; \der':= \s (\der ).$$
{\em Step 1}. $\s$ is a monomorphism.

$\noindent $

 {\em Step 2}. $\s (F) \subseteq F$,  where $F$ is
the only proper ideal of the algebra $\mI_1$. Therefore, there is
a
 commutative diagram   of algebra homomorphisms:
$$\xymatrix{\mI_1\ar[d]^\pi\ar[r]^{\s} & \mI_1\ar[d]^\pi\\
B_1\ar[r]^{\overline{\s} } &B_1}$$ where $B_1:= \mI_1/F\simeq
K[H][\der ,\der^{-1} ; \tau ]$, $\tau (H) = H+1$, is a simple
algebra, and so $\overline{\s}$ is an algebra monomorphism.

$\noindent $

{\em Step 3}. $H'=\l H+\mu +h$ for some elements $\l \in
K^*:=K\backslash \{ 0\}$, $\mu \in K$ and $h\in F$.

$\noindent $

{\em Step 4}. $H'=\frac{1}{n} H+\mu +f$, $\int' = \nu \int^n +f$
and $\der'= \nu^{-1} \der^n +g$ for some elements $\nu \in K^*$,
$n\geq 1$  and $h, f,g\in F$.

$\noindent $

{\em Step 5}. ${}^{\s}K[x]\simeq K[x]^n$, an isomorphism of
$\mI_1$-modules where $n$ is as in Step 4 and ${}_{\mI_1}K[x]:=
\mI_1/ \mI_1\der$, ${}^{\s}K[x]$ is the twisted $\mI_1$-module
$K[x]$ by the algebra endomorphism $\s$.

$\noindent $

{\em Step 6}. $n=1$, i.e.  ${}^{\s}K[x]\simeq K[x]$.

$\noindent $

{\em Step 7}. Up to the algebraic  torus action $\mT^1$
($\subseteq \Aut_{K-{\rm alg}}(\mI_1)$), $\nu=1$, i.e.
  $$H'= H+\mu +h, \;\; \int' =  \int +f, \;\;
 \der'=  \der +g.$$

{\em Step 8}. $\mu =0$.

$\noindent $

{\em Step 9}. $\s$ is an inner automorphism $\o_u$ of the algebra
$\mI_1$ for some unit  $u\in (1+F)^*$ of the algebra $\mI_1$.
$\Box $

$\noindent $

{\em Remark}. The algebra $B_1$ (see Step 3) is the left and right
localization of the Weyl algebra $A_1$ at the powers of the
element $\der$, i.e. the algebra $B_1$ is obtained from $A_1$ by
adding the {\em two-sided} inverse $\der^{-1}$ of the element
$\der$ (the algebra $B_1$ is also a {\em left} (but not right)
localization of the algebra $\mI_1$ at the powers of the element
$\der$, \cite{Bav-algintdif}, but in contrast to the Weyl algebra
$A_1$ the element $\der$ is not regular in $\mI_1$). An analogue
of the Conjecture/Problem of Dixmier fails for the algebra $B_1$:
for each natural number $n\geq 2$, the algebra monomorphism
$$ \s_n : B_1\ra B_1, \;\; H\mapsto \frac{1}{n}H, \;\; \der
\mapsto \der^n,$$ is obviously not an automorphism (use the
$\Z$-grading of the algebra $B_1= \bigoplus_{i\in \Z}K[H]\der^i$,
 $\der^i \alpha = \tau^i (\alpha ) \der^i$ for all $\alpha \in
 K[H]$ and $i\in \Z$). In view of  existence of this
 counterexample for the algebra $B_1$ it looks surprising that
 Theorem \ref{11Oct10} is true as the algebra $\mI_1$ is obtained
 from the Weyl algebra $A_1$ by adding a {\em right}, but not
 two-sided,  inverse of the element $\der$: $\der \int =1$ but
 $\int\der\neq 1$. Theorem \ref{11Oct10} can be seen  as a
 sign that the Conjecture/Problem of Dixmier is true.

$\noindent $

{\bf Conjecture}. {\em Each algebra endomorphism of $\mI_n$ is an
automorphism.}

$\noindent $

{\bf Ideas behind the proof of Theorem \ref{11Oct10}}. This is a
combination of old ideas/approach due to Dixmier \cite{Dix} of
using the eigenvalues of certain inner derivations (this was a key
moment in finding the group $\Aut_{K-{\rm alg}}(A_1)$ in
\cite{Dix} modulo many technicalities) and new ideas/approach of
using (i) the Fredholm operators and their indices based on the
fact that for the algebra $\mI_1$ the (Strong) Compact-Fredholm
Alternative holds \cite{Bav-intdifline} (which says that the
action of each polynomial integro-differential operator of $\mI_1$
on each simple $\mI_1$-module is either compact or Fredholm) and
(ii) the structure of the centralizers of elements of $\mI_1$
\cite{Bav-intdifline}.

$\noindent $

{\bf The Problem/Conjecture of Dixmier: recent progress}.
 In 1982, it was proved that a positive answer to the Problem/Conjecture of Dixmier
 for the  Weyl algebra $A_n$ implies the Jacobian Conjecture for
 the polynomial algebra $P_n$ in $n$ variables,
 see Bass, Connel and Wright
 \cite{BCW}. In 2005, it was proved independently by Tsuchimoto
 \cite{Tsuchi05} and Belov-Kanel and Kontsevich \cite{Bel-Kon05JCDP},
  see also \cite{JC-DP} for a short proof, that these two problems are equivalent.
 The Problem/Conjecture of Dixmier can be formulated as a question of
whether certain modules $\CM$ over the Weyl algebras are {\em
simple} \cite{Bav-Annals-2001} (recall that due to {\em Inequality
of Bernstein} \cite{Bernstein-1972} each simple module over the
Weyl algebra $A_n$ has the Gelfand-Kirillov dimension which is one
of the natural numbers $n, n+1, \ldots , 2n-1$;   Bernstein and
Lunts \cite{Bernstein-Lunts-1988}, \cite{Lunts-holonomic-1989}
showed that `generically' a simple $A_n$-module has the
Gelfand-Kirillov dimension $2n-1$). It is not obvious from the
outset that the modules $\CM$ are even finitely generated. In
2001, giving a positive answer to the Question of Rentschler on
the Weyl algebra  it was proved that the modules $\CM$ are
finitely generated and have the {\em smallest possible}
Gelfand-Kirillov dimension, i.e. $n$ (i.e. they are {\em
holonomic}) and as the result they have {\em finite length},
\cite{Bav-Annals-2001}. This means that the next step, as far as
the Jacobian Conjecture and the Problem/Conjecture of Dixmier are
concerned,  is either to prove the conjectures or to give a
counter-example.

One may wonder that for two different classes of algebras, the
polynomial algebras and the Weyl algebras, seemingly unrelated and
formulated in completely different ways conjectures, the Jacobian
Conjecture and the Conjecture of Dixmier, turned out to be
equivalent. It is obvious that there is   a phenomenon not yet
well understood. One may wonder that there are more algebras for
which one can formulate `similar' conjectures. Surprisingly, there
is a definite answer to this question:  in the class of all the
associative algebras conjecture like the two mentioned conjectures
makes sense {\em only} for the algebras $P_m\t A_n$ as was proved
in \cite{Bav-cpinv} (where $P_m$ is a polynomial algebra in $m$
variables; the two conjectures can be reformulated in terms of
locally nilpotent derivations that satisfy certain conditions, and
the algebras $P_m\t A_n$  are the only associative algebras that
have such derivations). This general conjecture for the algebras
$P_m\t A_n$  is true iff either the Jacobian Conjecture or the
Problem/Conjecture of Dixmier is true, see \cite{Bav-cpinv}.

$\noindent $

{\bf Meaning of the Problem/Conjecture of Dixmier and  the
Jacobian Conjecture, the groups of automorphisms}. The groups of
automorphisms of the polynomial algebra $P_n=P_1^{\t n}$, the Weyl
algebra $A_n=A_1^{\t n}$ and the algebra $\mI_n:=\mI_1^{\t n}
=K\langle x_1, \ldots , x_n, \frac{\der}{\der x_1}, \ldots
,\frac{\der}{\der x_n},  \int_1, \ldots , \int_n\rangle $ of
polynomial integro-differential operators are {\em huge}
infinite-dimensional algebraic groups. The groups of automorphisms
are known only for the  polynomial algebras when $n=1$ (trivial)
and $n=2$ (Jung  (1942) \cite{jung} and  Van der Kulk (1953)
\cite{kulk}); and for the Weyl algebra $A_1$ (Dixmier (1968)
\cite{Dix}) (in characteristic $p>0$, the group $\Aut_{K-{\rm
alg}}(A_1)$ was found by Makar-Limanov (1984)
\cite{Mak-LimBSMF84}, see also \cite{Bav-GMJ-09} for further
developments and another proof). In 2009, the group $\rG_n:=
\Aut_{K-{\rm alg}}(\mI_n)$ of automorphisms of the algebra $\mI_n$
 was  found for all $n\geq 1$, \cite{Bav-indifaut}:

\begin{eqnarray*}
\rG_n&=&S_n\ltimes \mT^n\ltimes  \Inn (\mI_n) \supseteq
 S_n\ltimes \mT^n \ltimes \underbrace{\GL_\infty (K)\ltimes\cdots \ltimes
\GL_\infty (K)}_{2^n-1 \;\; {\rm times}}, \\
\rG_1&\simeq & \mT^1 \ltimes \GL_\infty (K),
\end{eqnarray*}
 where  $S_n$
is the symmetric group, $\mT^n$ is the $n$-dimensional algebraic
torus,  $\Inn (\mI_n)$ is the group of inner automorphisms of
$\mI_n$ (which is huge). The ideas and approach in finding the
groups $\rG_n$ are completely different from that of Jung, Van der
Kulk and Dixmier: the Fredholm operators, ${\rm K}_1$-theory,
 indices.  On the other hand, when we look at the groups of
 automorphisms of the algebras $P_2$ , $A_1$ and $\mI_1$ (the
 only cases where we know explicit generators) we see that they
  have the `same nature': they are generated by affine
  automorphisms and `transvections.'

$\noindent $

The Jacobian Conjecture and the Problem/Conjecture of Dixmier (if
true) would give the {\em `defining relations'}  for the infinite
dimensional algebraic groups of automorphisms as infinite
dimensional varieties in the same way as the condition $\det =1$
defines the special linear (finite dimensional) algebraic group
$SL_n(K)$. Even true the conjectures would tell us nothing about
generators of the groups of automorphisms (i.e. about the
solutions of the defining relations, in the  same way and the
defining relation $\det =1$ tells nothing about generators for the
group $\SL_n(K)$).

$\noindent $

More obvious meaning of the Problem/Conjecture of Dixmier is that
the Weyl algebras $A_n$, which are simple {\em infinite
dimensional} algebras, behave like simple {\em finite dimensional}
algebras (each algebra endomorphism of a simple finite dimensional
algebra is, by a trivial reason,  an automorphism).  For a
polynomial algebra $P_n$ there are plenty algebra endomorphisms
that are not automorphisms. Recall that the {\em Jacobian
Conjecture} claims {\em that each algebra endomorphism $\s$ of the
polynomial algebra $P_n$ with the Jacobian ${\rm Jac} (\s ) :=\det
(\frac{\der \s ( x_i)}{\der x_j})\in K^*:=K\backslash \{ 0\}$ is
necessarily an automorphism}. The Jacobian condition is obviously
holds for all automorphisms of $P_n$ and the Jacobian Conjecture
implies that $\s$ is a {\em monomorphism}. So, the Jacobian
Conjecture (if true) means that each algebra monomorphism of $P_n$
which is as close as possible to be an automorphism {\em is}, in
fact, an automorphism.

$\noindent $

The paper is organized as follows. In Section \ref{AA1}, necessary
facts for the algebra $\mI_1$ are gathered which  are used later
in the paper. In Section \ref{PTH11}, the proof of Theorem
\ref{11Oct10} is given.


\section{The algebra $\mI_1$ }\label{AA1}

In this section, we collect necessary (mostly elementary)  facts
on the algebra $\mI_1$ from \cite{Bav-algintdif},
\cite{Bav-indifaut}, and \cite{Bav-intdifline} that are used later
in the paper.

The algebra $\mI_1$  is generated by the elements $\der $, $H:=
\der x$ and $\int$ (since $x=\int H$) that satisfy the defining
relations (Proposition 2.2, \cite{Bav-algintdif}):
\begin{equation}\label{I1rel}
\der \int = 1, \;\; [H, \int ] = \int, \;\; [H, \der ] =-\der ,
\;\; H(1-\int\der ) =(1-\int\der ) H = 1-\int\der ,
\end{equation}
 where
$[a,b]:=ab-ba$ is the {\em commutator} of elements $a$ and $b$.
 The elements of the algebra $\mI_1$,  
\begin{equation}\label{eijdef}
e_{ij}:=\int^i\der^j-\int^{i+1}\der^{j+1}, \;\; i,j\in \N ,
\end{equation}
satisfy the relations $e_{ij}e_{kl}=\d_{jk}e_{il}$ where $\d_{jk}$
is the Kronecker delta function and $\N := \{ 0, 1, \ldots \}$ is
the set of natural numbers. Notice that
$e_{ij}=\int^ie_{00}\der^j$. The matrices of the linear maps
$e_{ij}\in \End_K(K[x])$ with respect to the basis $\{ x^{[s]}:=
\frac{x^s}{s!}\}_{s\in \N}$ of the polynomial algebra $K[x]$  are
the elementary matrices, i.e.
$$ e_{ij}*x^{[s]}=\begin{cases}
x^{[i]}& \text{if }j=s,\\
0& \text{if }j\neq s.\\
\end{cases}$$
Let $E_{ij}\in \End_K(K[x])$ be the usual matrix units, i.e.
$E_{ij}*x^s= \d_{js}x^i$ for all $i,j,s\in \N$. Then
\begin{equation}\label{eijEij}
e_{ij}=\frac{j!}{i!}E_{ij},
\end{equation}
 $Ke_{ij}=KE_{ij}$, and
$F:=\bigoplus_{i,j\geq 0}Ke_{ij}= \bigoplus_{i,j\geq
0}KE_{ij}\simeq M_\infty (K)$, the algebra (without 1) of infinite
dimensional matrices. $F$ is the only proper ideal (i.e. $\neq 0,
\mI_1$) of the algebra $\mI_1$ \cite{Bav-algintdif}.  Using
induction on $i$ and the fact that $\int^je_{kk}\der^j=e_{k+j,
k+j}$, we can easily  prove that 
\begin{equation}\label{Iidi}
\int^i\der^i = 1-e_{00}-e_{11}-\cdots - e_{i-1,
i-1}=1-E_{00}-E_{11}-\cdots -E_{i-1,i-1}, \;\; i\geq 1.
\end{equation}

The monoid $1+F= 1+\bigoplus_{i,j\in \N} KE_{ij}=
1+\bigoplus_{i,j\in \N} Ke_{ij}$ admits the {\em determinant} map:
\begin{equation}\label{detFEij}
\det : 1+F \ra K, \;\; 1+\sum_{i,j=0}^d \l_{ij}E_{ij}\mapsto \det
(\sum_{i=0}^dE_{ii}+\sum_{i,j=0}^d \l_{ij}E_{ij}).
\end{equation}
By (\ref{eijEij}), this map can be defined as follows
\begin{equation}\label{1detFEij}
\det : 1+F \ra K, \;\; 1+\sum_{i,j=0}^d \l_{ij}e_{ij}\mapsto \det
(\sum_{i=0}^de_{ii}+\sum_{i,j=0}^d \l_{ij}e_{ij}).
\end{equation}
For all elements $a,b\in 1+F$, $\det (ab) = \det (a) \det (b)$ and
$ \det (1) = 1$. Therefore, an element $a\in 1+F$ is a unit iff
$\det (u)\neq 0$ (use the fact that $F$ is an ideal of $\mI_1$).

$\noindent $

{\bf $\Z$-grading on the algebra $\mI_1$ and the canonical form of
an integro-differential operator \cite{Bav-algintdif},
\cite{Bav-intdifline}}. The algebra $\mI_1=\bigoplus_{i\in \Z}
\mI_{1, i}$ is a $\Z$-graded algebra ($\mI_{1, i} \mI_{1,
j}\subseteq \mI_{1, i+j}$ for all $i,j\in \Z$) where
$$ \mI_{1, i} =\begin{cases}
D_1\int^i=\int^iD_1& \text{if } i>0,\\
D_1& \text{if }i=0,\\
\der^{|i|}D_1=D_1\der^{|i|}& \text{if }i<0,\\
\end{cases}
 $$
 the algebra $D_1:= K[H]\bigoplus \bigoplus_{i\in \N} Ke_{ii}$ is
a {\em commutative non-Noetherian} subalgebra of $\mI_1$, $
He_{ii} = e_{ii}H= (i+1)e_{ii}$  for $i\in \N $ (and so
$\bigoplus_{i\in \N} Ke_{ii}$ is the direct sum of non-zero ideals
$Ke_{ii}$ of the algebra $D_1$); $(\int^iD_1)_{D_1}\simeq D_1$,
$\int^id\mapsto d$; ${}_{D_1}(D_1\der^i) \simeq D_1$,
$d\der^i\mapsto d$,   for all $i\geq 0$ since $\der^i\int^i=1$.
 Notice that the maps $\cdot\int^i : D_1\ra D_1\int^i$, $d\mapsto
d\int^i$,  and $\der^i \cdot : D_1\ra \der^iD_1$, $d\mapsto
\der^id$, have the same kernel $\bigoplus_{j=0}^{i-1}Ke_{jj}$.

Each element $a$ of the algebra $\mI_1$ is the unique finite sum
\begin{equation}\label{acan}
a=\sum_{i>0} a_{-i}\der^i+a_0+\sum_{i>0}\int^ia_i +\sum_{i,j\in
\N} \l_{ij} e_{ij}
\end{equation}
where $a_k\in K[H]$ and $\l_{ij}\in K$. This is the {\em canonical
form} of the polynomial integro-differential operator
\cite{Bav-algintdif}.

$\noindent $

{\it Definition}. Let $a\in \mI_1$ be as in (\ref{acan}) and let
$a_F:=\sum \l_{ij}e_{ij}$. Suppose that $a_F\neq 0$ then
\begin{equation}\label{degFa}
\deg_F(a) :=\min \{ n\in \N \, | \, a_F\in \bigoplus_{i,j=0}^n
Ke_{ij}\}
\end{equation}
is called the $F$-{\em degree} of the element $a$;
$\deg_F(0):=-1$.

$\noindent $

Let $v_i:=\begin{cases}
\int^i& \text{if }i>0,\\
1& \text{if }i=0,\\
\der^{|i|}& \text{if }i<0.\\
\end{cases}$
Then $\mI_{1,i}=D_1v_i= v_iD_1$ and an element $a\in \mI_1$ is the
unique  finite  sum 
\begin{equation}\label{acan1}
a=\sum_{i\in \Z} b_iv_i +\sum_{i,j\in \N} \l_{ij} e_{ij}
\end{equation}
where $b_i\in K[H]$ and $\l_{ij}\in K$. So, the set $\{ H^j\der^i,
H^j, \int^iH^j, e_{st}\, | \, i\geq 1; j,s,t\geq 0\}$ is a
$K$-basis for the algebra $\mI_1$. The multiplication in the
algebra $\mI_1$ is given by the rule:
$$ \int H = (H-1) \int , \;\; H\der = \der (H-1), \;\; \int e_{ij}
= e_{i+1, j}, \;\; e_{ij}\int= e_{i,j-1}, \;\; \der e_{ij}=
e_{i-1, j}\;\; e_{ij} \der = \der e_{i, j+1}.$$
$$ He_{ii} = e_{ii}H= (i+1)e_{ii}, \;\; i\in \N, $$
where $e_{-1, j}:=0$ and $e_{i,-1}:=0$.

 $\noindent $

The factor algebra $B_1:= \mI_1/F$ is the simple Laurent skew
polynomial algebra $K[H][\der, \der^{-1}; \tau ]$ where the
automorphism $\tau \in \Aut_{K-{\rm alg}}(K[H])$ is defined by the
rule $\tau (H) = H+1$, \cite{Bav-algintdif}. Let
\begin{equation}\label{piI1B1}
\pi : \mI_1\ra B_1, \;\; a\mapsto \oa : a+F,
\end{equation}
be the canonical epimorphism.

$\noindent $

{\bf The groups of units $\mI_1^*$  and automorphisms
$\Aut_{K-{\rm alg}}(\mI_1)$ of the algebra $\mI_1$}. For a group
$G$, let $Z(G)$ denote its centre. Let $\mI_1^*$ be the group of
units of the algebra $\mI_1$. Since $F$ is an ideal of the algebra
$\mI_1$, the intersection $(1+F)^*:= \mI_1^* \cap (1+F)$ is a {\em
subgroup} of the group $\mI_1^*$. Moreover,
$$ (1+F)^* = \{ u\in 1+F \, | \, \det (u) \neq 0\}\simeq \GL_\infty (K).$$
The group $\Aut_{K-{\rm alg}}(\mI_1)$ of automorphisms of the
algebra $\mI_1$ contains the {\em algebraic torus}
$$ \mT^1:= \{ t_\l \, | \, \l \in K^*, \; t_\l (\int ) = \l \int,
\; t_\l (\der ) = \l^{-1} \der, \; t_\l (H) = H\} \simeq K^*, \;\;
t_\l \lra \l ,$$ and the group of inner automorphisms $\Inn
(\mI_1) = \{ \o_u : a\ra uau^{-1} \, | \, u\in \mI_1^*\}$ of the
algebra $\mI_1$.

\begin{theorem}\label{14Oct10}
\begin{enumerate}
\item {\rm (Theorem 4.5, \cite{Bav-algintdif})} $\mI_1^* =
K^*\times (1+F)^* \simeq K^*\times \GL_\infty (K)$ and $Z(\mI_1^*)
= K^*$. \item {\rm (Theorem 5.5.(1), \cite{Bav-indifaut})}
$\Aut_{K-{\rm alg}}(\mI_1)=\mT^1\ltimes  \Inn (\mI_1)$. \item {\rm
(Theorem 3.1.(2), \cite{Bav-indifaut})} The map $(1+F)^*\ra \Inn
(\mI_1)$, $u\mapsto \o_u$, is a group isomorphism.
\end{enumerate}
\end{theorem}



\section{Proof of Theorem \ref{11Oct10}}\label{PTH11}

This entire section is the proof of Theorem \ref{11Oct10}. We
follow the steps outlined in the Introduction.

Let $\s$ be an algebra endomorphism of $\mI_1$. We have to show
that $\s$ is an automorphism. The endomorphism $\s$ is uniquely
determined by its action on the generators $H$, $\int$ and $\der$
of the algebra $\mI_1$:
$$H':=\s (H), \;\; \int':= \s (\int ), \;\; \der':= \s (\der ).$$
{\em Step 1}. $\s$ {\em is a monomorphism.}

$\noindent $

Suppose that $\s$ is not a monomorphism, we seek a contradiction.
Then $\ker (\s )=F$ since $F$ is  the only proper (i.e. $\neq 0,
\mI_1$) ideal of the algebra $\mI_1$, \cite{Bav-algintdif}, and so
there is the algebra homomorphism
$$\overline{\s}:B_1:=\mI_1/F \ra \mI_1,\;\;
a+F\mapsto \s (a).$$  Since the algebra $B_1$ is a simple algebra,
$\overline{\s}$ is a monomorphism. The element $\der$ of the
algebra $B_1$ is an invertible  element and $\dim_K(K[\der
])=\infty$. Then $\overline{\s}=\s (\der )$ is an  invertible
element of the algebra $\mI_1$ and $\dim_K(K\langle \overline{\s}
(\der ) \rangle )= \dim_K(\overline{\s}(K[\der ] )=\dim_K(K[\der ]
)=\infty$ since $\overline{\s}$ is a monomorphism. This
contradicts the following lemma.

\begin{lemma}\label{a11Oct10}
For all units $u\in \mI_1^*$, $\dim_K(K\langle u \rangle ) <\infty
$.
\end{lemma}

{\it Proof}. The result follows from the equality $\mI_1^* = K^*
(1+F)^*$ (Theorem 4.5, \cite{Bav-algintdif}). $\Box $

$\noindent $

Therefore, $\s$ is a monomorphism.

$\noindent $

{\em Step 2}. $\s (F) \subseteq F$.

$\noindent $

\begin{lemma}\label{c11Oct10}
$K+F= \{ a\in \mI_1\, | \, \dim_K(K\langle a\rangle ) <\infty \}$.
\end{lemma}

{\it Proof}. The inclusion $\subseteq $ is obvious.  To show that
the inverse inclusion holds it suffices to prove  that,  for all
elements $a\not\in K+F$, $\dim_K(K\langle a \rangle )=\infty$, but
this is obvious since $ \oa:= a+F\in B_1\backslash K$ and
$\dim_K(K\langle \oa \rangle ) = \infty$. $\Box $

$\noindent $

By Lemma \ref{c11Oct10}, $\s (F) \subseteq K+F$. To prove that the
inclusion $\s (F) \subseteq F$ holds we have to show that $\s
(e_{ij}) \in F$ for all $i,j\in \N$. If $i=j$ then $e_{ii}^2=
e_{ii}$. If $\s (e_{ii})\not\in F$ then necessarily $\s
(e_{ii})\in \l_i+F$ for some $\l_i\in K^*$ such that $\l_i^2=
\l_i$, i.e. $\l_i= \pm 1$, we seek a contradiction. Since $\s
(K+F) \subseteq K+F$ and
\begin{eqnarray*}
\infty &=& \dim_K (\ker_{K+F}(\cdot e_{ii}))= \dim_K( \ker _{\s
(K+F)}(\cdot \s (e_{ii})))\;\;\;\;\;\; {\rm (by \; Step \; 1)}\\
&\leq & \dim_K( \ker_{K+F}(\cdot \s (e_{ii})))=\dim_K(
\ker_{K+F}(\cdot (\pm 1+f)))\\
&<&\infty ,
\end{eqnarray*}
 a contradiction. Then $\s (e_{ii})\in F$ for all $i\in \N$.

 For all $i\neq j$, $e_{ij}^2=0$, hence $\s (e_{ij})^2=0$, and so
 $\s (e_{ij})\in F$ since $\mI_1/ F$ is a domain. This proves that the
 inclusion $\s (F) \subseteq F$ holds. Therefore, there is a
 commutative diagram   of algebra homomorphisms:
$$\xymatrix{\mI_1\ar[d]^\pi\ar[r]^{\s} & \mI_1\ar[d]^\pi\\
B_1\ar[r]^{\overline{\s} } &B_1}$$
 where $\overline{\s}(a+F) = \s (a) +F$ for all $a\in \mI_1$; $\pi
 : \mI_1\ra B_1=\mI_1/ F$, $a\mapsto a+F$; and so $\overline{\s}$ is
 an algebra monomorphism since $B_1$ is a simple
algebra.

$\noindent $

{\em Step 3}. $H'=\l H+\mu +h$ {\em for some elements $\l \in
K^*:=K\backslash \{ 0\}$, $\mu \in K$ and $h\in F$ where $F$ is
the only proper ideal of the algebra $\mI_1$.}

$\noindent $

For an element $a\in \mI_1$, let $\Cen_{\mI_1}(a) = \{ b\in \mI_1,
\, | \, ab=ba \}$ be its {\em centralizer} in the algebra $\mI_1$,
and $\Cen_F(a) := F\cap \Cen_{\mI_1}(a)$.

\begin{proposition}\label{Ca1Jun10}
{\rm \cite{Bav-intdifline}} Let $a\in \mI_1$. Then
$\dim_K(\Cen_F(a)) =\infty$ iff $a\in K[H]+F$.
\end{proposition}


By Proposition \ref{Ca1Jun10},
 $H'\in K[H]+F$, i.e. $H'=\alpha +h$ for unique elements
 $\alpha \in K[H]$ and $h\in F$ since $K[H]\cap F=0$ (see (\ref{acan1})). Since, for
 each element $\th \in \{ H,  \int , \der \}$,
 $$ \infty = \dim_K(K[\th ]) = \dim_K(\s (K[\th ] )= \dim_K(K[\s
 (\th )])\;\; {\rm and}\;\; \dim_K(K\langle \l +f\rangle )<\infty, $$
 for all elements $\l \in K$ and $f\in F$,
 we must have
\begin{equation}\label{KF1}
\alpha \in K[H]\backslash K \;\; {\rm and}\;\; \int',
 \der'\not\in K+F.
\end{equation}
Using (\ref{I1rel}) and the direct sum decomposition
 $$\mI_1= \bigoplus_{i\geq 1} D_1\der^i\bigoplus D_1\bigoplus
 \bigoplus_{i\geq 1} \int^i D_1,$$
  we see that the set of eigenvalues of the  inner derivation $\ad (H) : \mI_1\ra \mI_1$, $\mapsto
 [H,a]:= Ha-aH$, of the algebra $\mI_1$ is $\Ev (\ad (H))=\Z$,
 and, for each eigenvalue $i\in \Z$,
 $$\ker_{\mI_1}(\ad (H) - i) = \begin{cases}
\int^iD_1& \text{if }i\geq 1,\\
D_1& \text{if }i=0,\\
D_1\der^{|i|}& \text{if }i\geq -1.\\
\end{cases}
$$
Since $\s$ is a monomorphism (by Step 1), $\Ev (\ad (H'))\supseteq
\Z$. By (\ref{KF1}), $\pi (\int'^i) \neq 0$ and $\pi (\der'^i)\neq
0$ for all $i\in \N$ where $\pi$ is defined in (\ref{piI1B1}).
Since, by (\ref{I1rel}),
$$ [ \pi (H'),\pi (\int'^i)]=i\pi (\int'^i)\;\; {\rm and}\;\;  [ \pi (H'),\pi (\der'^i)]=-i\pi
(\der'^i), \;\; {\rm for \; all\;}\; i\geq 1,$$ we see that $\Ev
(\pi (H') , B_1)\supseteq \Z$. By (\ref{KF1}), $\pi (H') =\alpha
\in K[H]\backslash K$.

By Lemma \ref{b11Oct10}, 
\begin{equation}\label{KF2}
\alpha = \l H +\mu
\end{equation}
for some $\l \in K^*$ and $\mu \in K$.
\begin{lemma}\label{b11Oct10}
Let $a\in K[H]\backslash K$. Then $\Ev (\ad (a) , B_1) \neq 0$ iff
$a= \l H+\mu$ where $\l \in K^*$ and $\mu \in K$.
\end{lemma}

{\it Proof}. ($\Leftarrow $) Obvious: $[\l H+\mu , \der ] = -\l
\der$.

($\Rightarrow $) It suffices to show that if $\deg_H(a)>1$ then
$\Ev (\ad  (a) , B_1)=0$. The algebra $B_1= \bigoplus_{i\in \Z}
K[H]\der^i$ is a $\Z$-graded algebra where $K[H]\der^i$ is the
$i$'th graded component of the algebra $B_1$. The element $a\in
K[H]$ is a homogeneous element of the algebra $B_1$. Therefore,
for each eigenvalue $\nu \in \Ev (\ad (a) , B_1)$,
$$ \ker_{B_1}(\ad (a) - \nu ) = \bigoplus_{i\in \Z} (\ker_{B_1}(\ad
(a) - \nu ) \bigcap K[H]\der^i).$$ Suppose that $\nu \neq 0$, then
$[a, \beta \der^i] =\nu  \beta \der^i$ for some elements $0\neq
\beta \in K[H]$ and $i\in \Z$,  necessarily $i\neq 0$ since $\nu
\neq 0$. The equality can be written as $(a-\tau^i (a))\beta
\der^i = \nu \beta \der^i$, and so $ a-\tau^i ( a)=\nu$ since
$B_1$ is a domain. Since $\deg_H(a-\tau^i (a)) = \deg_H( a) -
1\geq 1$, this is impossible. Therefore, $\Ev (\ad (a) , B_1) =0$.
$\Box $

$\noindent $

{\em Step 4}. $H'=\frac{1}{n} H+\mu +h$, $\int' = \nu \int^n +h$
{\em and $\der'= \nu^{-1} \der^n +g$ for some elements $\nu \in
K^*$, $n\geq 1$  and} $h, f, g\in F$.

$\noindent $

By Step 3, $\Ev ( \pi (H') =\l H +\mu , B_1) = \l \Ev ( H, B_1) =
\l \Z$ and, for each element $i\in \Z$,
$$ \ker_{B_1} (\ad (\pi ( H')) - i\l ) = B_1 \der^{-i}.$$
Applying the algebra homomorphism $\pi \s$ to the relations $[H,
\int ] =\int$, $[H, \der ] = -\der$ and $\der \int =1$ yields the
equalities
$$[\pi (H'), \pi (\int') ] =\pi (\int'), \;\; [\pi (H'), \pi (\der') ]
= -\pi (\der') , \;\; \pi (\der') \pi ( \int') =1.$$ By
(\ref{KF1}), $\pi (\int') \neq 0$ and $\pi (\der') \neq 0$.
Therefore, by Lemma \ref{b11Oct10}, there are two options
\begin{eqnarray*}
 &(i)&  \pi (H') = \frac{1}{n}H+\mu, \;\; \pi (\int') = \nu \der^{-n},
 \;\; \pi (\der') = \nu^{-1} \der^n; \\
 &(ii)& \pi (H') = -\frac{1}{n}H+\mu, \;\; \pi (\int') = \nu^{-1} \der^n,
 \;\; \pi (\der') = \nu \der^{-n};
\end{eqnarray*}
for some natural number $n\geq 1$ and $\nu \in K^*$ since
$$B_1=\bigoplus_{i\in \Z} K[H]\der^i, \;\; \Ev (\ad (H), B_1) =
\Z, \;\; \ker (\ad (H)-i)=K\der^{-i}, \;\; i\in \Z.$$ Therefore,
there are elements $h,f,g\in F$ such that
\begin{eqnarray*}
 &(i)&  H' = \frac{1}{n}H+\mu +h, \;\; \int' = \nu \int^n+f,
 \;\; \der' = \nu^{-1} \der^n+g; \\
 &(ii)& H' = -\frac{1}{n}H+\mu +h, \;\; \int' = \nu^{-1} \der^n+g,
 \;\; \der' = \nu \int^n+f.
\end{eqnarray*}
We are going to show that the case (ii) is not possible. For we
need some results.

$\noindent $

Since $ \der' \int'=1$, the map $\der'\cdot : K[x]\ra K[x]$,
$p\mapsto \der'*p$, is a {\em surjection}, and so 
\begin{equation}\label{dkid}
\dim_K(\ker_{K[x]}(\der'\cdot )) = \ind_{K[x]}(\der'\cdot )
\end{equation}
where $\ind_{K[x]}(\v ):= \dim_K(\ker_{K[x]}(\v )) -
\dim_K(\coker_{K[x]}(\v ))$ is the {\em index} of a linear map $\v
\in \End_K(K[x])$ provided the  kernel and cokernel of the map $\v
$ are finite dimensional.

\begin{theorem}\label{CB30May10}
{\rm (\cite{Bav-intdifline})} Let $a\in \mI_1$, $M$ be a nonzero
$\mI_1$-module of finite length and $a_M: M\ra M$, $m\mapsto am$.
Then $\dim_K(\ker (a_M))<\infty $ iff $\dim_K(\coker (a_M))<\infty
$ iff $a\not\in F$.
\end{theorem}
\begin{lemma}\label{D3Oct10}
{\rm (\cite{Bav-intdifline})} Let $a\in \mI_1\backslash F$ and
$f\in F$. Then $\ind_M(a+f) = \ind_M(a)$ for all left or right
$\mI_1$-modules $M$ of finite length where $\ind_M(a):=
\dim_K(\ker(a_M))-\dim_K(\coker (a_M))$.
\end{lemma}


{\em Step 5}. ${}^{\s}K[x]\simeq K[x]^n$, {\em  an isomorphism of
$\mI_1$-modules where $n$ is as in Step 4}.

$\noindent $

Here ${}_{\mI_1}K[x]:= \mI_1/ \mI_1\der$ is a faithful
$\mI_1$-module (since $\mI_1\subseteq \End_K(K[x])$), and the
action of an element $a\in \mI_1$ on a polynomial $p\in K[x]$ is
denoted by $a*p$.  ${}^{\s}K[x]$ is the twisted by the algebra
endomorphism $\s$  $\mI_1$-module $K[x]$: as vector spaces ${}^\s
K[x]=K[x]$ but the action of the algebra $\mI_1$ on ${}^\s K[x]$
is given by the rule, $a\cdot p := \s (a) *p$ for all elements
$a\in \mI_1$ and $p\in K[x]$. The $\mI_1$-module $K[x]$ is a
simple (since $A_1\subseteq \mI_1$ and the $A_1$-module $K[x]$ is
simple),  and $\der'\not\in F$, by (\ref{KF1}). By (\ref{dkid})
and  Lemma \ref{D3Oct10},
$$ \ker_{K[x]} (\der'\cdot ) = \ind_{K[x]}(\der'\cdot ) =
\ind_{K[x]}((\nu^{-1} \der^n + g)\cdot ) = \ind_{K[x]} (\der^n )
=\ker_{K[x]}(\der^n \cdot ) = n$$ since $\ker_{K[x]}(\der^n \cdot
) = \bigoplus_{i=0}^{n-1} Kx^i$. Recall that ${}_{\mI_1}K[x]\simeq
\mI_1/ \mI_1\der$ is a simple $\mI_1$-module such that
\begin{equation}\label{Kxx}
K[x]=\bigcup_{i\geq 1} \ker_{K[x]}(\der^i\cdot ),\;
\;\ker_{K[x]}(\der^i\cdot )=\bigoplus_{j=0}^{i-1} Kx^j, \;\;
\ker_{K[x]}(\der\cdot )=K.
\end{equation}
Similarly, for a natural number $n\geq 1$, the direct sum $K[x]^n$
of $n$ copies of the simple $\mI_1$-module $K[x]$ is a semi-simple
$\mI_1$-module of finite length $n$, $K[x]^n = \bigcup_{i\geq 1}
\ker_{K[x]^n} (\der^i\cdot )$ and $\dim_K( \ker_{K[x]^n} ( \der
\cdot )=n\dim_{K[x]}(\der \cdot) = n$. It follows that the
$\mI_1$-module epimorphism
$$ \v : K[x]^n \ra V:= \mI_1\cdot \ker_{ {}^\s K[x]}(\der \cdot ) =
\s ( \mI_1) * \ker_{K[x]}( \der'\cdot )$$ (where
${}_{\mI_1}V\subseteq {}_{\mI_1}({}^\s K[x])$) given by the rule
$$ \v : (1, 0, \ldots , 0)\mapsto v_1, \ldots , (0, \ldots ,
0,1)\mapsto v_n,$$ is an {\em isomorphism} where $(1, 0, \ldots ,
0), \ldots , (0\ldots , 0,1)$ is the standard free $K[x]$-basis
for the $\mI_1$-module $K[x]^n$ and $v_1, \ldots , v_n$ is a
$K$-basis for the vector space $ \ker_{K[x]}(\der'\cdot )$
(otherwise, ${}_{\mI_1}V\simeq K[x]^m$ for some $m<n$, and so
$n=\dim_K(\ker_V(\der'\cdot))=  \dim_K(\ker_{K[x]^m} ( \der'\cdot
) )=m$, a contradiction.)

Fix $s\in \N$ such that $$s>\max \{ n,d, \deg_F(h), \deg_F(f),
\deg_F(g)\}$$ where $\ker_{K[x]}(\der'\cdot ) \subseteq K[x]_{\leq
d}:=\bigoplus_{i=0}^dKx^i$ for some number $d\in \N$. Then, for
all integers $i\geq s$, by Step 4,
\begin{eqnarray*}
 H'*x^{[i]}&=& (\frac{1}{n}(i+1)+\mu ) x^{[i]},  \\
\int'*x^{[i]}&=&  \nu x^{[i+n]}, \\
\der'*x^{[i]}&=& \nu^{-1} x^{[i-n]},\\
\end{eqnarray*}
where $x^{[j]}:=0$ for all integers $j<0$. For each integer $i\in
\N$, let $K[x]_{\leq i}:=\bigoplus_{j=0}^i Kx^j$. Then
$K[x]=\bigcup_{i\in \N} K[x]_{\leq i}$.  Consider the ascending
chain of vector spaces in $K[x]$:
$$ V_0:=K[x]_{\leq s}\subset V_1:= K[x]_{\leq s+n}\subset \cdots
\subset V_t:=K[x]_{\leq s+nt}\subset \cdots ,$$  $\dim_K(V_t) =
1+s+nt$. Then, for all $t\in \N$,
\begin{eqnarray*}
 H'*V_t&\subseteq & V_t,  \\
\int'*V_t&\subseteq & V_{t+1}, \\
\der'*V_t&\subseteq &V_t.\\
\end{eqnarray*}
Since $\ker_{K[x]}(\der'\cdot ) \subseteq V_0$ and $K[x]_{\leq
t}=\sum_{i=0}^tK\int^i*1$, we see that $ \v (K[x]^n_{\leq
t})\subseteq \sum_{i=0}^t K\int'^i*V_0 \subseteq V_{\leq t}$, and
so
$$\dim_K(V_{\leq t})-\dim_K(\v (K[x]^n_{\leq t})) = 1+s+nt-n(t+1)
= 1+s-n={\rm const}.$$ This means that the factor $\mI_1$-module
${}^\s K[x]/ \im (\v )$ is {\em finite dimensional}. Therefore,
$$ \im (\v ) = {}^\s K[x],$$
since the only finite dimensional $\mI_1$-module is the zero one
(the algebra $\mI_1$ contains the simple infinite dimensional
algebra $A_1$, and the only finite dimensional $A_1$-module is the
zero one), i.e. the $\mI_1$-modules $K[x]^n$ and ${}^\s K[x]$ are
isomorphic via $\v$.

$\noindent $

{\em Step 6}. $n=1$, i.e.  ${}^{\s}K[x]\simeq K[x]$.

$\noindent $

By Step 5,  ${}^{\s}K[x]\simeq K[x]^n$. Notice that
$K[x]=\bigoplus_{i\in \N} Kx^i$ and $Kx^i = \ker_{K[x]}(H-i-1)$,
i.e. the linear map $H\cdot : K[x]\ra K[x]$, $p\mapsto H*p$, is
semi-simple. Therefore, the map $H\cdot : K[x]^n\ra K[x]^n$,
$p\mapsto H\cdot p$, is semi-simple and each of its eigenvalues
has {\em multiplicity} (i.e. the dimension of the corresponding
eigenspace) $n$. Since  ${}_{\mI_1}({}^{\s}K[x])\simeq K[x]^n$
(Step 5), the linear map $H'\cdot : K[x]\ra K[x]$, $p\mapsto
H'*p$, is semi-simple and each its eigenvalue has multiplicity
$n$.  Since
\begin{eqnarray*}
H'*x^{[i]}&=& (\frac{1}{n}(i+1)+\mu ) x^{[i]},\;\;\;\;\; i\geq s,   \\
H'*V_0&\subseteq & V_0, \;\;\;\;  \dim_K(V_0)<\infty, \\
\end{eqnarray*}
we must have
$$n=1,$$ (since the eigenvalues $\{ \frac{1}{n}(i+1)+\mu \, | \, i\geq s\}$ of
the linear map $H'\cdot$ acting in $K[x]$
 are all distinct) and so
$$H'=H+\mu +h, \;\; \int'= \nu \int +f, \;\; \der'= \nu^{-1} \der
+g.$$ Up to the algebraic  torus action $\mT^1$ ($\subseteq
\Aut_{K_{\rm alg}}(\mI_1)$), we may assume that $\nu=1$, i.e.

$\noindent $

{\em Step 7}. $H'= H+\mu +h$, $\int' =  \int +f$ and $ \der'=
 \der +g$.

$\noindent $

{\em Step 8}. $\mu =0$.

$\noindent $

For the $\mI_1$-module $K[x]$ and for all natural numbers $i\geq
1$, 
\begin{equation}\label{idimH1}
K[x]_{\leq i-1}=\ker_{K[x]}(\der^i\cdot )=
\bigoplus_{j=1}^{i}\ker_{K[x]}(H-j)=\der K[\der ]*
\ker_{K[x]}(H-(i+1))
\end{equation}
 and 
\begin{equation}\label{idimH}
i=\dim_K(\der K[\der ]* \ker_{K[x]}(H-(i+1)))\;\; {\rm for \;
all}\;\; i\in \Ev (H\cdot , K[x])=\{ 1,2, \ldots \}.
\end{equation}
 Since the vector space
$U:= V_0\bigoplus Kx^{[s+1]}= K[x]_{\leq s+1}$ is
$\der'$-invariant, $\der'*V_0\subseteq V_0$,
$\der'*x^{[s+1]}=x^{[s]}\in V_0$, 
\begin{equation}\label{HpS1}
 H'*x^{[s+1]}= (s+1+1+\mu )x^{[s+1]}\;\; {\rm and } \;\;
{}_{\mI_1}({}^\s K[x])\simeq K[x]
\end{equation}
 we must have, by (\ref{idimH1}),
$$ V_0= \der'K[\der']*\ker_{K[x]}(H'-(s+2+\mu))=\der'K[\der']*x^{[s+1]}.$$
 By (\ref{idimH}) and since
${}_{\mI_1}({}^\s K[x])\simeq K[x]$,
 $$ (s+2+\mu ) -1= \dim_K(\der' K[\der']*\ker_{K[x]}(H'-(s+2+\mu
 ))) = \dim_K(V_0)= \dim_K(K[x]_{\leq s}=s+1.$$Therefore, $\mu =0$.

$\noindent $

{\em Step 9}. $\s$ {\em is an inner automorphism $\o_u$ of the
algebra $\mI_1$ for some unit  $u\in (1+F)^*$ of the algebra}
$\mI_1$.

$\noindent $

Notice that
\begin{eqnarray*}
  K[x]&=&\bigoplus_{i\geq 1} \ker_{K[x]}(H'-i), \;\; \Ev (H'\cdot, K[x]) = \{ 1,2, \ldots , \}, \\
  \dim_K (\ker_{K[x]}(H'-i ))&=&1\;\; {\rm for \; all } \; i\in
 \Ev (H'\cdot, K[x]), \\
  \int'*\ker_{K[x]}(H'-i)&=&\ker_{K[x]}(H'-(i+1)) \;\; {\rm
  and}\;\; \\
 \der'*\ker_{K[x]}(H'-i)&=&\ker_{K[x]}(H'-(i-1))\;\; {\rm for \; all } \; i\in
 \Ev (H'\cdot, K[x]).
\end{eqnarray*}
Since $$K[x]=V_0\bigoplus
(x^{s+1})=\bigoplus_{i=0}^sKx'^{[i]}\bigoplus Kx^{[s+1]}\bigoplus
Kx^{[s+2]}\bigoplus+\cdots $$ where $(x^{s+1})=K[x]x^{s+1}$,
$x'^{[i]}:= \der'^{(s+1-i)}*x^{[s+1]}$ and ${}_{\mI_1}({}^\s
K[x])\simeq K[x]$, we see that (by (\ref{HpS1}))
$$ \ker_{K[x]}(H'-i-1)=
\begin{cases}
Kx'^{[i]}& \text{if }i=0,1, \ldots , s,\\
Kx^{[i]}& \text{if }i>s.\\
\end{cases}$$
Let $x'^{[i]}:=x^{[i]}$ for all $i>s$. Then
$$\der'*x'^{[i]}=
\begin{cases}
x'^{[i-1]}& \text{if }i>0,\\
0& \text{if }i=0.\\
\end{cases}$$
Then necessarily,
$$\int'*x'^{[i]}=x'^{[i+1]}, \;\; i\geq 0, $$
using the facts that $\der'\int'=1$,   $\int'*\ker_{K[x]}(H'-i) =
\ker_{K[x]}(H'-i-1)$ and $\ker_{K[x]}(H'-i)=Kx'^{[i-1]}$ for all
$i\geq 1$. The $K$-linear map
$$ u: K[x]\ra K[x], \;\; x^{[i]}\mapsto x'^{[i]}, $$
is an $\mI_1$-module isomorphism $u: K[x]\ra {}^\s K[x]$ since
$$ ua*x^{[i]}= a'*(ux^{[i]})$$
for all elements $a\in \{ H,\int , \der \}$ and $i\in \N$, i.e.
$ua = \s (a) u$, and so $\s (a) = uau^{-1}= \o _u (a)$ for all
elements $a\in \{ H, \int , \der \}$. Notice that $u\in (1+F)^*$,
i.e. $\s = \o_u\in \Inn (\mI_1)$. $\Box $

$\noindent $

$${\bf Acknowledgements}$$

$\noindent $

I would like to thank J. Dixmier for inspiring discussions on his
problems and  comments during the Solstice Conference in June
2010, Paris. I would like to thank P. Adjamagbo, J.-Y. Charbonnel
and  R. Rentschler for interesting discussions.

Department of Pure Mathematics

University of Sheffield

Hicks Building

Sheffield S3 7RH

UK

email: v.bavula@sheffield.ac.uk

\end{document}